\documentclass[12pt]{article}
\usepackage{amsmath,amssymb,amsbsy,amsfonts,amsthm,latexsym,amsopn,amstext,
                                          amsopn,amstext,amsxtra,euscript,amscd}
\begin{document}

\newfont{\teneufm}{eufm10}
\newfont{\seveneufm}{eufm7}
\newfont{\fiveeufm}{eufm5}
%
%
\newfam\eufmfam
                          \textfont\eufmfam=\teneufm
\scriptfont\eufmfam=\seveneufm
                          \scriptscriptfont\eufmfam=\fiveeufm

%
%
\def\frak#1{{\fam\eufmfam\relax#1}}
%


\def\bbbr{{\rm I\!R}} 
\def\bbbc{{\rm I\!C}} 
\def\bbbm{{\rm I\!M}}
\def\bbbn{{\rm I\!N}} 
\def\bbbf{{\rm I\!F}}
\def\bbbh{{\rm I\!H}}
\def\bbbk{{\rm I\!K}}
\def\bbbl{{\rm I\!L}}
\def\bbbp{{\rm I\!P}}
\newcommand{\lcm}{{\rm lcm}}
\def\bbbone{{\mathchoice {\rm 1\mskip-4mu l} {\rm 1\mskip-4mu l}
{\rm 1\mskip-4.5mu l} {\rm 1\mskip-5mu l}}}
\def\bbbc{{\mathchoice {\setbox0=\hbox{$\displaystyle\rm C$}\hbox{\hbox
to0pt{\kern0.4\wd0\vrule height0.9\ht0\hss}\box0}}
{\setbox0=\hbox{$\textstyle\rm C$}\hbox{\hbox
to0pt{\kern0.4\wd0\vrule height0.9\ht0\hss}\box0}}
{\setbox0=\hbox{$\scriptstyle\rm C$}\hbox{\hbox
to0pt{\kern0.4\wd0\vrule height0.9\ht0\hss}\box0}}
{\setbox0=\hbox{$\scriptscriptstyle\rm C$}\hbox{\hbox
to0pt{\kern0.4\wd0\vrule height0.9\ht0\hss}\box0}}}}
\def\bbbq{{\mathchoice {\setbox0=\hbox{$\displaystyle\rm
Q$}\hbox{\raise
0.15\ht0\hbox to0pt{\kern0.4\wd0\vrule height0.8\ht0\hss}\box0}}
{\setbox0=\hbox{$\textstyle\rm Q$}\hbox{\raise
0.15\ht0\hbox to0pt{\kern0.4\wd0\vrule height0.8\ht0\hss}\box0}}
{\setbox0=\hbox{$\scriptstyle\rm Q$}\hbox{\raise
0.15\ht0\hbox to0pt{\kern0.4\wd0\vrule height0.7\ht0\hss}\box0}}
{\setbox0=\hbox{$\scriptscriptstyle\rm Q$}\hbox{\raise
0.15\ht0\hbox to0pt{\kern0.4\wd0\vrule height0.7\ht0\hss}\box0}}}}
\def\bbbt{{\mathchoice {\setbox0=\hbox{$\displaystyle\rm
T$}\hbox{\hbox to0pt{\kern0.3\wd0\vrule height0.9\ht0\hss}\box0}}
{\setbox0=\hbox{$\textstyle\rm T$}\hbox{\hbox
to0pt{\kern0.3\wd0\vrule height0.9\ht0\hss}\box0}}
{\setbox0=\hbox{$\scriptstyle\rm T$}\hbox{\hbox
to0pt{\kern0.3\wd0\vrule height0.9\ht0\hss}\box0}}
{\setbox0=\hbox{$\scriptscriptstyle\rm T$}\hbox{\hbox
to0pt{\kern0.3\wd0\vrule height0.9\ht0\hss}\box0}}}}
\def\bbbs{{\mathchoice
{\setbox0=\hbox{$\displaystyle     \rm S$}\hbox{\raise0.5\ht0\hbox
to0pt{\kern0.35\wd0\vrule height0.45\ht0\hss}\hbox
to0pt{\kern0.55\wd0\vrule height0.5\ht0\hss}\box0}}
{\setbox0=\hbox{$\textstyle        \rm S$}\hbox{\raise0.5\ht0\hbox
to0pt{\kern0.35\wd0\vrule height0.45\ht0\hss}\hbox
to0pt{\kern0.55\wd0\vrule height0.5\ht0\hss}\box0}}
{\setbox0=\hbox{$\scriptstyle      \rm S$}\hbox{\raise0.5\ht0\hbox
to0pt{\kern0.35\wd0\vrule height0.45\ht0\hss}\raise0.05\ht0\hbox
to0pt{\kern0.5\wd0\vrule height0.45\ht0\hss}\box0}}
{\setbox0=\hbox{$\scriptscriptstyle\rm S$}\hbox{\raise0.5\ht0\hbox
to0pt{\kern0.4\wd0\vrule height0.45\ht0\hss}\raise0.05\ht0\hbox
to0pt{\kern0.55\wd0\vrule height0.45\ht0\hss}\box0}}}}
\def\bbbz{{\mathchoice {\hbox{$\sf\textstyle Z\kern-0.4em Z$}}
{\hbox{$\sf\textstyle Z\kern-0.4em Z$}}
{\hbox{$\sf\scriptstyle Z\kern-0.3em Z$}}
{\hbox{$\sf\scriptscriptstyle Z\kern-0.2em Z$}}}}
\def\ts{\thinspace}

\newtheorem{theorem}{Theorem}
\newtheorem{lemma}[theorem]{Lemma}
\newtheorem{claim}[theorem]{Claim}
\newtheorem{cor}[theorem]{Corollary}
\newtheorem{prop}[theorem]{Proposition}
\newtheorem{definition}{Definition}
\newtheorem{question}[theorem]{Open Question}

\def\squareforqed{\hbox{\rlap{$\sqcap$}$\sqcup$}}
\def\qed{\ifmmode\squareforqed\else{\unskip\nobreak\hfil
\penalty50\hskip1em\null\nobreak\hfil\squareforqed
\parfillskip=0pt\finalhyphendemerits=0\endgraf}\fi}

\def\cA{{\mathcal A}}
\def\cB{{\mathcal B}}
\def\cC{{\mathcal C}}
\def\cD{{\mathcal D}}
\def\cE{{\mathcal E}}
\def\cF{{\mathcal F}}
\def\cG{{\mathcal G}}
\def\cH{{\mathcal H}}
\def\cI{{\mathcal I}}
\def\cJ{{\mathcal J}}
\def\cK{{\mathcal K}}
\def\cL{{\mathcal L}}
\def\cM{{\mathcal M}}
\def\cN{{\mathcal N}}
\def\cO{{\mathcal O}}
\def\cP{{\mathcal P}}
\def\cQ{{\mathcal Q}}
\def\cR{{\mathcal R}}
\def\cS{{\mathcal S}}
\def\cT{{\mathcal T}}
\def\cU{{\mathcal U}}
\def\cV{{\mathcal V}}
\def\cW{{\mathcal W}}
\def\cX{{\mathcal X}}
\def\cY{{\mathcal Y}}
\def\cZ{{\mathcal Z}}

\newcommand{\comm}[1]{\marginpar{%
\vskip-\baselineskip 
\raggedright\footnotesize
\itshape\hrule\smallskip#1\par\smallskip\hrule}}





\hyphenation{re-pub-lished}

\def\ord{{\mathrm{ord}}}
\def\Nm{{\mathrm{Nm}}}
\renewcommand{\vec}[1]{\mathbf{#1}}

\def \F{{\bbbf}}
\def \L{{\bbbl}}
\def \K{{\bbbk}}
\def \Z{{\bbbz}}
\def \N{{\bbbn}}
\def \Q{{\bbbq}}
\def\E{{\mathbf E}}
\def\bH{{\mathbf H}}
\def\G{{\mathcal G}}
\def\O{{\mathcal O}}
\def\cS{{\mathcal S}}
\def \R{{\bbbr}}
\def\Fp{\F_p}
\def \fp{\Fp^*}
\def\\{\cr}
\def\({\left(}
\def\){\right)}
\def\fl#1{\left\lfloor#1\right\rfloor}
\def\rf#1{\left\lceil#1\right\rceil}

\newcommand{\li}{\operatorname{li}}

\def\Zm{\Z_m}
\def\Zt{\Z_t}
\def\Zp{\Z_p}
\def\Um{\cU_m}
\def\Ut{\cU_t}
\def\Up{\cU_p}

\def\ep{{\mathbf{e}}_p}

\def \Prob{{\mathrm {}}}

\def\LC{{\cL}_{C,\cF}(Q)}
\def\LCn{{\cL}_{C,\cF}(nG)}
\def\Mrs{\cM_{r,s}\(\F_p\)}

\def\Fbar{\overline{\F}_q}
\def\Fn{\F_{q^n}}
\def\En{\E(\Fn)}

\def\kk{\mathbf{k}}
\def\hh{\mathbf{h}}
\def\eps{\varepsilon}

\def\mand{\qquad \mbox{and} \qquad}

\def\MOV{{\bf{MOV}}}


\title{On Curves over Finite Fields with Jacobians of Small Exponent}

\author{
{\sc Kevin Ford} \\
{Department of Mathematics, 1409 West Green Street} \\
{University of Illinois at Urbana-Champaign} \\
{Urbana, IL 61801, USA} \\
{\tt ford@math.uiuc.edu} \\
\and
{\sc Igor Shparlinski} \\
{Department of Computing, Macquarie University} \\
{Sydney, NSW 2109, Australia} \\
{\tt igor@ics.mq.edu.au}
}

\date{\today}

\maketitle

\begin{abstract}
We show that finite fields over which there is a
curve of a given genus $g\ge 1$ with its Jacobian having a small exponent,
are very rare. This extends a recent result of W.~Duke in the case
$g=1$.  We also show when $g=1$ or $g=2$,
our lower bounds on the exponent, valid for almost
all finite fields $\F_q$ and all curves over $\F_q$, are best possible.
\end{abstract}



\section{Introduction}

Let $\cJ_\cC(\F_q)$ denote the Jacobian of a curve $\cC$
defined over a finite field $\F_q$ of $q$ elements.
We denote by $\ell_q(C)$ the exponent
of $\cJ_\cC(\F_q)$ (that  is, $\ell_q(\cC)$ is the largest order of
elements of the Abelian group $\cJ_\cC(\F_q)$) and by $g$ the
genus of $\cC$.
We start with recalling two well know facts.
\begin{itemize}
\item The Weil bound implies that
      \begin{equation}
\label{eq:Jac Card}
(q^{1/2}-1)^{2g} \le \# \cJ_\cC(\F_q) \le (q^{1/2}+1)^{2g},
\end{equation}
see Corollary~5.70, Theorem~5.76 and Corollary~5.80 of~\cite{ACDFLNV}.
In particular, for fixed $g$,
$$
\# \cJ_\cC(\F_q) = q^g + O_g(q^{g-1/2}).
$$
\item The Jacobian  $ \cJ_\cC(\F_q)$ is an Abelian group with
at most $2g$ generators, that is,
for some positive integers $m_1, \ldots, m_{2g}$ we have
      \begin{equation}
\label{eq:Jac Struct}
      \cJ_\cC(\F_q) \cong \Z/m_1\Z \times  \ldots \times\Z/m_{2g}\Z,
\quad
\text{where}
\quad m_1\mid  \ldots \mid m_{2g},
\end{equation}
(in particular  $m_1=\ldots=m_j=1$ if the rank of
$ \cJ_\cC(\F_q)$ is $2g - j$)
and also
\begin{equation}\label{mgq1}
m_i | (q-1) \qquad (1\le i\le g),
\end{equation}
see Proposition~5.78 of~\cite{ACDFLNV}.
\end{itemize}

Thus we see $\ell_q(\cC) = m_{2g}$ where $ m_{2g}$ is defined by the
representation~\eqref{eq:Jac Struct}, which together with~\eqref{eq:Jac Card}
implies the following trivial bound
      \begin{equation}
\label{eq:Exp Trivial}
\ell_q(\cC) \ge (\# \cJ_\cC(\F_q) )^{1/2g} \ge q^{1/2}-1.
\end{equation}

For elliptic curves $\cC = \cE$
over finite fields the exponent $\ell_q(\cE)$ has been
studied in a number of works, see~\cite{Duke,LuMcKSh,LuSh,Schoof,Shp},
with a variety of  results, each of them indicating that in
a ``typical case''  $\ell_q(\cE)$  tends to be substantially
larger than the bound~\eqref{eq:Exp Trivial} guarantees.
However for general curves the behavior of $\ell_q(\cC)$
has not been studied. Let $\pi(x)$ denote the number
of primes $p \le x$.
      W.~Duke~\cite[footnote on page~691]{Duke}, among
other results, has proved that  for a sufficiently large $x$ and all
but $o(\pi(x))$ of prime powers $q \le x$, the bound
\begin{equation}
\label{eq:Duke's bound}
\ell_q(\cE) \ge  q^{3/4}/\log q
\end{equation}
holds for all elliptic curves $\cE$ defined over $\F_q$ (the paper~\cite{Duke}
considers only primes, but including all prime powers in the
statement is trivial of course).

We provide a generalization and some improvement
of~\eqref{eq:Duke's bound} for curves of arbitrary genus.

\begin{theorem}
\label{thm:Exp Large} Fix $g\ge 1$ and let $\eps(x)$ be a
positive, decreasing function of $x$ with $\eps(x)\to 0$ as $x\to\infty$.
For all but $o(\pi(x))$ of the prime powers $q \le x$, the bound
$$
\ell_q(\cC) \ge  q^{3/4 + \eps(q)}
$$
holds for all curves $\cC$ of genus $g$ defined over $\F_q$.
\end{theorem}

The method of proof of~\eqref{eq:Duke's bound},
used in~\cite{Duke}, is somewhat specific to  elliptic
curves, so here we use a slightly different approach
to counting fields $\F_q$ that may contain a ``bad'' curve.

We show that Theorem \ref{thm:Exp Large} is best possible
for $g=1$ and $g=2$.  In particular,
the bound~\eqref{eq:Duke's bound} of W.~Duke~\cite{Duke} is quite sharp.

%
%

\begin{theorem}
\label{thm:Duke's bound} 
For any fixed $\varepsilon > 0$ there exists $\alpha > 0$
such that for  sufficiently large $x$, there are at least
$\alpha \pi(x)$ primes $q \le x$ such that for some
nonsupersingular elliptic curve $\cE$ and some
nonsupersingular curve $\cC$ of genus $g=2$
defined over $\F_q$, the bounds
$$
\ell_q(\cE) \le  q^{3/4 + \varepsilon} \mand \ell_q(\cC) \le  q^{3/4
+ \varepsilon}
$$
hold.
\end{theorem}

The proof is based on a special case of a certain lower
bound on the number of
shifted primes $p-1$ having a divisor in a given interval.
In full generality this bound is given in Theorem~7 of~\cite{Ford}.
Such results have been applied  to study the order of a
given integer $a> 1$ modulo almost all primes $p$,
see~\cite{ErdMur,IndlTim,Papp}, and now they have turned out to be
useful for studying exponents of  Jacobians.
This argument also immediately implies the following result which
applies to all curves over $\F_q$ of all possible genera.

\begin{theorem}
\label{thm:Exp Large Any g} Let $\eps(x)$ be a
positive, decreasing function of $x$ with $\eps(x)\to 0$ as $x\to\infty$.
For all but $o(\pi(x))$ of the prime powers $q \le x$, the bound
$$
\ell_q(\cC) \ge  q^{1/2 + \eps(q)}
$$
holds for all  curves $\cC$ of arbitrary genus defined over $\F_q$.
\end{theorem}

Throughout the paper, the implied constants in the symbols `$O$',
`$\ll$' and `$\gg$' do not depend on any parameter unless indicated
by a subscript, that is, $O_g$, $\ll_g$ or $\gg_g$
(we recall that the notations $U = O(V)$,
$U\ll V$, and $V \gg U$ are all equivalent to the assertion
that the inequality $|U|\le cV$ holds for some constant $c>0$).

\bigskip

{\bf Acknowledgments.} The authors would like to thank Bjorn Poonen for
a number of very helpful discussions.
During the preparation of this paper,
K.~F.\ was supported by NSF grants DMS-0301083 and DMS-0555367
and I.~S.\ was supported in part by ARC grant DP0556431.

\section{Preliminaries}

We have already mentioned that our results are based on
some estimates from~\cite{Ford} on shifted primes having a
divisor in a given interval. Here we give a brief guide to
these estimates.

As in~\cite{Ford} we use $H(x,y,z)$ to denote the number of
positive integers $n \le x$ having a divisor $d$ with $y < d \le z$.
Theorem~1 of~\cite{Ford} gives the right order of magnitude of
$H(x,y,z)$ in the full range of parameters. However for our purposes
we need only the estimate
\begin{equation}
\label{eq:Hxyz}
H(x,y,z) \ll x u^\delta (\log (2/u))^{-3/2}
\end{equation}
where
$$
\delta   = 1 -\frac{1 + \log \log 2}{\log 2}
= 0.086071\ldots
$$
and $u$ is defined by the equation $y^{1 + u} = z$, which holds
uniformly in the range $2y \le z \le y^2$, $3\le y\le \sqrt{x}$.

Furthermore, we need the upper bound on $H(x,y,z)$ only as
tool of estimating $H(x,y,z, \cP_\lambda)$ which is the
number of primes $p \le x$ such that $p+\lambda$
has a divisor $d$ with $y < d \le z$.
Theorem~6 of~\cite{Ford} gives the upper bound
\begin{equation}
\label{eq:HxyzP}
H(x,y,z, \cP_\lambda) \ll \frac{H(x,y,z)}{\log x}
\end{equation}
which holds for every fixed non-zero integer $\lambda$ in the
range $z \ge y + (\log y)^{2/3}$ and $3\le y \le \sqrt{x}$,
which is much wider than is
necessary for the purposes of this paper.

We also need Theorem~7 of~\cite{Ford} which gives a lower bound
on $H(x,y,z, \cP_\lambda) $ in a certain range of $x,y,z$.
However, since   its proof is quite short, we give
an independent derivation in Section~\ref{sec:Thm 2}.

\section{Proof of Theorem~\ref{thm:Exp Large}}

The number of prime powers $q=p^a\le x$ with $a\ge 2$ is $O(x^{1/2})$.
Thus, it suffices to show that for
all but $o(x/\log x)$ of the primes $q$ with $x/2 < q \le x$,
     the bound
$$
\ell_q(\cC) \ge  q^{3/4 + \eps(q)}
$$
holds for all  curves $\cC$ of genus $g$ defined over $\F_q$.

For a $(2g-1)-$tuple $\kk=(k_1,\ldots,k_{2g-1})$ of
positive integers,
we consider the set $\cQ_{\kk}$  of primes
$x/2 \le q \le x$ for which there exists a curve
$\cC$ of genus $g\ge 1$ over $\F_q$ such that  $m_1 = k_1$,
$m_i = m_{i-1} k_i$,
where  $m_i$ is as in~\eqref{eq:Jac Struct} and~\eqref{mgq1},
$i =1, \ldots, 2g-1$.  In particular, if such a curve $\cC$ exists,
then
\begin{equation}
\label{eq: q-1 congr}
q-1 \equiv 0 \pmod{k_1 \ldots k_g}.
\end{equation}


Since
$$
k_1^{2g} k_2^{2g-1} \ldots k_{2g-1}^2 | \# \cJ_\cC(\F_q),
$$
we see by~\eqref{eq:Jac Card} that there are at most
\begin{equation}
\label{eq:Set N'}
U_{\kk} = \frac{(x^{1/2}+1)^{2g}} {k_1^{2g} k_2^{2g-1} \ldots k_{2g-1}^2}
\end{equation}
possibilities for the cardinality $N = \# \cJ_\cC(\F_q)$.

For each of such values $N$, we see by~\eqref{eq:Jac Card} that
$$
N^{1/g}-2N^{1/2g}+1 \le q \le N^{1/g} + 2N^{1/2g} + 1.
$$
Recalling~\eqref{eq: q-1 congr} we deduce that for each possible cardinality
$N$ the prime powers $q$ may take at most
\begin{equation}
\label{eq:Set q'}
V_{\kk} =   \frac{5(x^{1/2}+1)}{k_1 k_2\ldots k_g} + 1
\end{equation}
values.
Therefore, combining~\eqref{eq:Set N'} and~\eqref{eq:Set q'}, we derive
\begin{equation}\label{eq:Qki}
\# \cQ_{\kk} \le  U_{\kk}  V_{\kk} \le
\frac{5(x^{1/2}+1)^{2g+1}} {k_1^{2g+1} k_2^{2g} \ldots k_g^{g+2}k_{g+1}^{g}
      \ldots k_{2g-1}^2}  +
\frac{(x^{1/2}+1)^{2g}} {k_1^{2g} k_2^{2g-1} \ldots k_{2g-1}^2}.
\end{equation}
When $g=1$, we interpret the right side as $5(x^{1/2}+1)^3 k_1^{-3}+
(x^{1/2}+1)^2 k_1^{-2}$.

For any curve $\cC$ of genus
$g\ge 1$ over $\F_q$ and any positive integer $s \le 2g-1$,
we have
\begin{equation}\label{eq:m2g lower}
\ell_q(\cC) =  m_{2g} \ge \(\frac{\# \cJ_\cC(\F_q)} {m_1\ldots
m_s}\)^{1/(2g-s)}
\ge \(\frac{(q^{1/2}-1)^{2g}}{k_1^s k_2^{s-1} \ldots k_s}\)^{1/(2g-s)}.
\end{equation}
In fact, we  only need~\eqref{eq:m2g lower} for $s=g$ and $s=2g-1$.

Suppose without loss of generality that $\eps(x) \ge (\log x)^{-1/2}$
and  write $\eta=\eps(x/2)$.  Assume $x$ is large, in
particular so large that
$$
\eta<\frac{1}{100g}.
$$
Let $I$ be the interval $(x^{1/4-3\eta}, x^{1/4+3\eta}]$.
Let $\cK$ denote the set of $\kk$ satisfying
\begin{align}
k_1\cdots k_g &\not\in I, \label{eq:k1} \\
k_1^g k_2^{g-1} \ldots k_g &\ge x^{g/4-2g \eta}, \label{eq:k2} \\
k_1^{2g-1} k_2^{2g-2} \ldots k_{2g-1} &\ge
x^{g-3/4-2\eta}. \label{eq:k3}
\end{align}

Partition the primes $q\in (x/2,x]$ into three sets: $\cT_1$ is the set
      of such primes for which $q-1$ has a divisor in $I$,
$\cT_2$ is the set of such
      primes lying in a set $Q_{\kk}$ with $\kk \in \cK$, and $\cT_3$ is the
      set of remaining primes.  By Theorems~1 and~6 of~\cite{Ford},
that is, by a combination of~\eqref{eq:Hxyz} and~\eqref{eq:HxyzP},
   we
      have
\begin{equation}\label{eq:T1}
\# \cT_1 \ll \frac{x}{\log x} \eta^{\delta} (\log 1/\eta)^{-3/2}
\end{equation}

Now consider $q\in \cT_2$.
By \eqref{eq:k2},
$$
k_1 \cdots k_g \ge (k_1^g k_2^{g-1} \cdots k_g)^{1/g} \ge x^{1/4-2\eta},
$$
hence $k_1\cdots k_g > x^{1/4+3\eta}$ by \eqref{eq:k1}.
Combined with \eqref{eq:Qki}, \eqref{eq:k3}, and the inequality
$k_i\le (x^{1/2}+1)^{2g}$ for each $i$, we obtain
\begin{align*}
\# \cT_2 &\le \sum_{\kk\in\cK} \# \cQ_{\kk} \\
&\le \( \frac{5(x^{1/2}+1)^{2g+1}}{x^{g-1/2+\eta}}
   + \frac{(x^{1/2}+1)^{2g}}{x^{g-3/4-2\eta}}\) \sum_{\kk\in\cK}
\frac{1}{k_1 \cdots k_{2g-1}} \\
&\le  \( \frac{5(x^{1/2}+1)^{2g+1}}{x^{g-1/2+\eta}}
   + \frac{(x^{1/2}+1)^{2g}}{x^{g-3/4-2\eta}}\) (2g\log(x^{1/2}+1)+1)^{2g-1} \\
&\ll_g (\log x)^{2g-1} ( x^{1-\eta} + x^{3/4+2\eta} ) \\
&\ll_g x^{1-\eta/2}.
\end{align*}

Together with~\eqref{eq:T1}, we see that all but $o(x/\log x)$ primes
$q\in (x/2,x]$ lie in $\cT_3$.  For $q\in \cT_3$,
the condition~\eqref{eq:k1} holds,
      thus either~\eqref{eq:k2} is false or~\eqref{eq:k3} is false.
In either case, the bound~\eqref{eq:m2g lower} implies that $\ell_q(\cC) \gg_g
x^{3/4+2\eta}$, and hence for large $x$
$$
\ell_q(\cC) \ge q^{3/4+\eps(q)}
$$
for any curve  $\cC$ of genus $g$ defined over $\F_q$.  \qed

%
%
\section{Proof of Theorem~\ref{thm:Duke's bound}}
\label{sec:Thm 2}
%
%

We start with the case $g=1$.

Without loss of generality we can assume that $\varepsilon < 1/20$.
Put
$$
y=x^{1/4-\varepsilon} \mand z=x^{1/4 - \varepsilon/2}.
$$
Since $y> x^{1/5}$, an integer $k \le x$
      can have at most 4 prime
factors $p$ with $y<p\le z$.
Hence, the set  $\cP$ of primes $x/\log x \le q \le x$ such that
$q-1$ has a prime
divisor $p$ with $y<p\le z$, is of cardinally least
$$
\# \cP \ge \frac{1}{4}  \sum_{\substack{y<p\le z\\p~\mathrm{prime}}}
\pi(x;p,1) +
O\(\frac{x}{(\log x)^2}\),
$$
where, as usual, $\pi(x;k,a)$ is the number of primes $q\le x$
with $q \equiv a \pmod k$.

By the Bombieri-Vinogradov theorem
(see, for example,
Section~28 of~\cite{Dav}),
$$
\sum_{\substack{y<p\le z\\p~\mathrm{prime}}}
\left| \pi(x;p,1) - \frac{1}{p-1}\pi(x)\right|  \ll \frac{x}{(\log x)^2} .
$$
Therefore
$$
\# \cP \ge \frac{1}{4}  \pi(x) \sum_{\substack{y<p\le
z\\p~\mathrm{prime}}} \frac{1}{p-1} + O\(\frac{x}{(\log x)^2}\)
=\frac{1}{4}  \pi(x) \sum_{\substack{y<p\le
z\\p~\mathrm{prime}}} \frac{1}{p} + O\(\frac{x}{(\log x)^2}\).
$$
By the Mertens  theorem (see Theorem~4.1
of Chapter~1 in~\cite{Prach}),
$$
\sum_{\substack{y<p\le
z\\p~\mathrm{prime}}} \frac{1}{p} = \log \log z - \log \log y + o(1)
= \log \frac{1 - 2 \varepsilon}{1 - 4 \varepsilon} + o(1),
$$
thus for large $x$ we have $\# \cP \ge  \alpha \pi(x)$ for a positive $\alpha$
depending on $\eps$.  This result is a special case of Theorem~7
of~\cite{Ford}, but we include the proof because it is short.

For a sufficiently large $x$ and for any $q \in \cP$, there are at least
$2q^{1/2} z^{-2} - 1 \ge q^\varepsilon$  integers $k \in [q+1- 2q^{1/2}, q]$
with $p^2 | k$ for some prime $p|q-1$ with $y < p \le z$.
For any such $k$, by~\cite{Ruck,TsfVlad,Voloch} one can always find an
elliptic curve $\cE$
over $\F_q$ with $\cE(\F_q) = k$ of $\F_q$-rational
points and the exponent $\ell_q(\cE) = k/p\le q/y
\le q^{3/4 + \varepsilon}$. This concludes the proof in the case $g=1$.

For $g=2$, Proposition~5.4 in  Section~5 of  Chapter~X
of~\cite{Silv} implies that the cardinalities of elliptic curves $\cE$ over
$\F_q$ with $j$-invariant $j(\cE)=0, 1728$ take $O(1)$ values.
Therefore we can choose $k$ and an elliptic curve $\cE$
over $\F_q$ of exponent $\ell_q(\cE) \le q^{3/4 + \varepsilon}$
as in the above with the
additional  condition  $j(\cE) \ne 0, 1728$.
By Corollary~6 of~\cite{HLP} we see that there is a curve $\cC$ of
genus $g=2$ such that the Jacobian  $J_\cC(\F_q)$ is isogenous to
$\cE(\F_q)\times \cE(\F_q)$. Moreover, there
exists an isogeny from  $\cE(\F_q)\times \cE(\F_q)$
to $J_\cC(\F_q)$,  whose kernel (over an algebraic closure of $\F_q$) is
isomorphic to $\Z/2\Z \times \Z/2\Z$.
So  $\ell_q(\cC) \ge  \ell_q(\cE)/2$, which concludes the proof for $g=2$.
\qed

\section{Proof of Theorem~\ref{thm:Exp Large Any g}}


The desired bound follows immediately from
     Theorems~1 and~6 of~\cite{Ford},
that is, from~\eqref{eq:Hxyz} and~\eqref{eq:HxyzP},
and the congruence $q-1 \equiv 0 \pmod {m_g}$,
where  $m_i$, $ i =1, \ldots, 2g$, are as in~\eqref{eq:Jac Struct}.
Again without loss of generality assume that
      $\eps(x) \ge (\log x)^{-1/2}$.
For $\eta=2 \eps(x/2)$, similarly to~\eqref{eq:T1},  we see that the
set $\cR$ of
primes
$q\le x$ such that $q-1$ has   a divisor $m \in [x^{1/2 - 2\eta}, x^{1/2 +
2\eta}]$, is of cardinality $\# \cR = o(x/\log x)$.
Consider a prime $q\in (2x^{1-\eta},x]$ which does not lie in $\cR$,
and any curve $\cC$ of genus $g$ over $\F_q$.  If $m_g>x^{1/2+\eta}$ then
$$
\ell_q(\cC) = m_{2g} \ge m_g > q^{1/2+\eps(q)}.
$$
Otherwise, by~\eqref{mgq1}, $m_g \le x^{1/2-2\eta}$ and
by~\eqref{eq:Jac Card}
we obtain
\begin{eqnarray*}
\ell_q(\cC) & \ge & \(\frac{\# \cJ_\cC(\F_q)}{m_1 \cdots m_g}\)^{1/g} \ge
    \(\frac{ (q^{1/2} -1)^{2g}}{m_1 \cdots m_g}\)^{1/g} \\
& \ge & \( \frac{x^{g-g\eta}}{x^{g/2 - 2g\eta}}\)^{1/g} \ge
x^{1/2+\eta} > q^{1/2+\eps(q)}
\end{eqnarray*}
for large $x$.
\qed

\section{Remarks}

It is interesting to note that
using~\eqref{eq:m2g lower} for other values of
$s$ (besides $s=g$ and $s=2g-1$ as in the proof of
Theorem~\ref{thm:Exp Large})
and thus corresponding  sets $\cK$, does not lead to any improvements.

\medskip

\noindent
{\bf Open Question.}
{\it
Is the exponent in Theorem~\ref{thm:Exp Large}
sharp for arbitrary $g\ge 3$, as it is for $g=1, 2$?}
\medskip

Unfortunately the lack of knowledge about the distribution
of possible cardinalities of Jacobians of curves of genus $g \ge 2$
prevents are from deriving an analogue of Theorem~\ref{thm:Duke's bound}
for $g \ge 2$.

%

\end{document}